\documentstyle[a4,amsthm,twoside,amssymb,12pt]{article}
\begin{document}
\newcommand{\Kdim}{\mathrm{K.dim}}
\newcommand{\lann}{\mathrm{l.ann}}
\newcommand{\rann}{\mathrm{r.ann}}
\renewcommand{\max}{\mathrm{max}}
\newcommand{\Ker}{\mathrm{Ker}\;}

\newtheorem{thm}{Theorem}
\newtheorem{prop}[thm]{Proposition}
\newtheorem{lem}[thm]{Lemma}
\newtheorem{cor}[thm]{Corollary}
\theoremstyle{definition}
\newtheorem*{defn}{Definition}
\newtheorem{exmp}{Example}
\theoremstyle{remark}
\newtheorem*{rem}{Remark}
 \newcommand{\ZZ}{{\mathbb Z}}
 \newcommand{\NN}{{\mathbb N}}

 \renewenvironment{proof}{\par\noindent{\bf Proof.}}{$\square$\par\bigskip}

\title{ MODULES WHOSE SMALL SUBMODULES HAVE KRULL DIMENSION}
\author{Christian Lomp\\
        Mathematisches Institut \\
        Heinrich-Heine Universit\"at D\"usseldorf\\
        40225 D\"usseldorf, Germany}
\date{ }
\maketitle

\begin{abstract}
The main aim of this paper is to show that an $AB5^{\ast}$ module
whose small submodules have Krull dimension has a radical having Krull
dimension. The proof uses the notion of dual Goldie dimension.
\end{abstract}

\section{Introduction}

Let $R$ be an associative ring with unit and let $M$ be a left unital
$R$-module. Denote the socle of $M$, the intersection of all essential
submodules of $M$, by $Soc(M)$. A well-known theorem by Goodearl (see
\cite[Proposition 3.6]{goodearl} or \cite[Proposition 4]{alkhazzi})
asserts that $M/Soc(M)$ is noetherian if and only if every factor module
$M/N$ with $N$ essential in $M$ is noetherian. This can easily be
extended to show that $M/Soc(M)$ has Krull dimension if and only if
$M/N$ has Krull dimension for every essential submodule $N$ of $M$ (see
\cite[Proposition 2]{puczylowski}).
Denote the radical of $M$, the sum of all small submodules of $M$, by $Rad(M)$.
Dual to Goodearl's result Al-Khazzi and Smith proved
that $Rad(M)$ is artinian if and only if every small submodule of $M$ is artinian
(see \cite[Theorem 5]{alkhazzi}). They asked in \cite{alkhazzi}:
If every small submodule has finite uniform dimension (Goldie
dimension). Does $Rad(M)$ have finite uniform dimension ?\\
Puczy{\l}owski answered this question in the negative and showed that there
exists a $\ZZ$-module $M$ such that every
small submodule is noetherian and hence has Krull dimension but $Rad(M)$
does not have Krull dimension (see \cite[Example]{puczylowski}).

Since we wish to dualize Goodearl's Theorem it is natural to
ask if the Al-Khazzi-Smith Theorem can be extended for arbitrary Krull
dimension to modules which satisfy property $AB5^{\ast}$.

Theorem \ref{main_theorem} is our main theorem and shows that for a module
having $AB5^{\ast}$ the following implication holds:
If every small submodules has finite hollow dimension (dual Goldie
dimension) then every submodule of $Rad(M)$ has finite hollow dimension.
This can be seen as dual to \cite[Lemma 5.14]{dung}: If $M/N$ has finite
uniform dimension for every essential submodule $N$ of $M$, then every
factor module of $M/Soc(M)$ has finite uniform dimension.

\section{Definitions}
For the definition of Krull dimension we refer to Chapter 6 in
\cite{dung}. A module $M$ is said to be {\it uniform} if $M\neq0$ and
every non-zero submodule is essential in $M$. $M$ is said to have
{\it finite uniform dimension}
(or {\it finite Goldie dimension})
if there is a monomorphism from a finite direct sum
of proper uniform submodules of $M$ to $M$ such that the image is essential
in $M$. It is well known that this is equivalent to the property that
$M$ has no infinite independent family of non-zero submodules and that
there is a maximal finite independent family of uniform submodules (see
\cite[Theorem 5.9]{dung}).
We denote the cardinality of this family by $udim(M)$ and
call $udim(M)$ the {\it uniform dimension} of $M$.

A module $M$ is said to be {\it hollow} if $M\neq0$ and every proper
submodule is small in $M$. Hollow modules were introduced by Fleury in
\cite{fleury}.
% More on hollow modules can be found in \cite{wisbauer}.
$M$ is said to have {\it finite hollow dimension} if there
is an epimorphism with a small kernel from $M$ to a finite direct sum of
non-zero hollow factor modules. It can be shown, that in this case
there is a number $n$ such that $M$ does not allow an epimorphism to a
direct sum with more than $n$ summands. We denote this by $hdim(M)=n$
and call $hdim(M)$ the {\it hollow dimension} of $M$. For any submodule
$N$ of $M$ we have $hdim(M/N) \leq hdim(M)$.

\begin{defn} Let $M$ be an $R$-module and $\{ N_\lambda \}_\Lambda$ a
family of proper submodules of $M$. $\{ N_\lambda \}_\Lambda$ is called
{\it coindependent} (see \cite{takeuchi}) if for every $\lambda \in \Lambda$ and finite subset
$J \subseteq \Lambda \setminus \{ \lambda \}$
$$ N_\lambda + \bigcap_{j \in J} N_j = M$$
holds (convention: if $J$ is empty, then set $\bigcap_J N_j := M$).
\end{defn}
It can be shown, that a module $M$ has finite hollow dimension if and only if
every coindependent family of submodules is finite (see \cite[Corollary
13]{gr} ).\\
For more information on dual Goldie dimension we refer to \cite{gr},
\cite{reiter}, \cite{takeuchi} and \cite{var1}.

\begin{defn}
 An $R$-module $M$ has property {\it $AB5^*$} if for every submodule
 $N$ and inverse systems $\{ M_i \}_{i \in I} $ of submodules of $M$ the
 following holds:
 \begin{center}
   $N + \bigcap_{i \in I}M_i = \bigcap_{i\in I}(N+M_i)$
 \end{center}
\end{defn}

Examples of modules with $AB5^{\ast}$ are artinian or
linearly compact modules. Herbera and Shamsuddin proved the following result:
\begin{lem} \label{hs} (\cite[Lemma 6]{herbera} or \cite{brodskii})\\
For a module $M$ with property $AB5^{\ast}$ the following statements are
equivalent:
\begin{enumerate}
\item[(a)] Every factor module of $M$ has finite uniform dimension.
\item[(b)] Every submodule of $M$ has finite hollow dimension.
\end{enumerate}
\end{lem}

It is easy to see that implication $(b) \Rightarrow (a)$ always holds
(see \cite[Proposition 12]{var2}) and that $(a) \Rightarrow (b)$ is
false in general (for example $M=_\ZZ\ZZ$).

\begin{defn}
Let $M$ be an $R$-module and $\{ N_\lambda \}_\Lambda$ a family of
proper submodules. Then $\{ N_\lambda \}_\Lambda$ is called
{\it completely coindependent} if for every $\lambda \in \Lambda$:
$$ N_\lambda + \bigcap_{\mu \neq \lambda} N_\mu = M$$
holds.
\end{defn}

A completely coindependent family is coindependent, but the converse is
not true in general (for example $\{ p\ZZ \}$ in $_\ZZ\ZZ$ where $p$ runs through all
prime numbers). Considering Herbera and Shamsuddin's proof of Lemma \ref{hs}
we get:
\begin{lem} \label{lemma1}
Every coindependent family of submodules of a module with property
$AB5^{\ast}$ is completely coindependent.
\end{lem}

\section{Modules whose small submodules have Krull dimension.}
In this section we will prove our main theorem. First we prove:

\begin{lem} \label{lemma3}
Let $M$ be an $R$-module, $\{ N_\lambda \}_\Lambda$ a
completely coindependent family of proper submodules of $M$ and
$|\Lambda| \geq 2$. Assume that for every
$\lambda \in \Lambda$ there exists a submodule $L_\lambda$ such that
$N_\lambda \subsetneqq L_\lambda$. Let
$L:=\bigcap_{\lambda \in \Lambda} L_\lambda$ and $N:=\bigcap_{\lambda
\in \Lambda} N_\lambda$.
Then $\{ (N_\lambda \cap L)/N \}_\Lambda$ forms a completely
coindependent family of proper submodules of $L/N$.
\end{lem}

\begin{proof} Let $\lambda \in \Lambda$. Then
$ N_\lambda + L =  L_\lambda \cap ( N_\lambda + \bigcap_{\mu \neq
\lambda} L_\mu) = L_\lambda \cap M =L_\lambda. $
Since $N_\lambda \neq L_\lambda$ we have $N_\lambda \cap L \subsetneqq L$.
Moreover
$(N_\lambda \cap L) + \bigcap_{\mu \neq \lambda} (N_\mu \cap L) = L$ is
straightforward.
Thus $\{ N_\lambda \cap L \}_\Lambda$ forms a completely
coindependent family of proper submodules of $L$.
Hence $N \subsetneqq N_\lambda \cap L$ for all $\lambda \in \Lambda$
and $\{ (N_\lambda \cap L)/N \}_\Lambda$
is a completely coindependent family of proper submodules of $L/N$.
\end{proof}

The next definition dualizes the notion of an essential extension.
\begin{defn}
Let $N \subseteq L \subseteq M$ be submodules of $M$. We say that $L$ {\it
lies above} $N$ (in M) if $L/N \ll M/N$.
Note that $L$ lies above $N$ if and only if
$N+K=M$ holds whenever $L+K=M$ holds for a submodule $K$ of $M$.
\end{defn}

\begin{lem} \label{lemma4}
 Let $M$ be an $R$-module with $AB5^{\ast}$,
 $\{ L_\lambda \}_\Lambda$ a coindependent family of submodules such that
 for each $\lambda \in \Lambda$ there exists a submodule $N_\lambda
 \subseteq L_\lambda$ such that $L_\lambda$ lies above $N_\lambda$ in
 $M$. Then $\bigcap_\Lambda L_\lambda$ lies above $\bigcap_\Lambda N_\lambda$ in
 $M$.
\end{lem}

\begin{proof}
Let $\Omega$ denote the set of all
finite subsets of $\Lambda$. Define for all $J \in \Omega$
$A_J:= \bigcap_{j \in J} L_j \mbox{ and } B_J:=\bigcap_{j \in J} N_J.$
By induction on the cardinality of $J$ it is easy to show, that
$A_J$ lies above $B_J$ for all $J\in \Omega$ (see \cite[Proposition
1.6]{takeuchi}). Since $\{A_J \}_{J \in \Omega}$ and $\{ B_J \}_{J \in \Omega}$ are inverse
systems, we get for $K\subset M$:
\begin{eqnarray*}
 M & = & K + \bigcap_{\lambda \in \Lambda} L_\lambda
     = K + \bigcap_{J \in \Omega} A_J
     = \bigcap_{J \in \Omega} (K + A_J)
     = \bigcap_{J \in \Omega} (K + B_J) \\
   & = & K + \bigcap_{J \in \Omega} B_J
     =  K + \bigcap_{\lambda \in \Lambda} N_\lambda.
\end{eqnarray*}
\end{proof}

\begin{defn}
Let $M$ be an $R$-module and $N,L$ submodules of $M$.
Then $N$ is called a {\it supplement} of $L$ in $M$
if $N$ is minimal with respect to $N+L=M$.
Note that $N$ is a supplement of $L$ in $M$ if and only if
$N+L=M$ and $N\cap L \ll N$ (see \cite[Chapter 41]{wisbauer}).
A module is called {\it amply supplemented} if whenever $N+L=M$ holds for two
submodules of $M$, then $N$ contains a supplement of $L$ in $M$.
Any module with $AB5^{\ast}$ is amply supplemented
(see \cite[47.9]{wisbauer}).
As a generalization of a supplement we say that
$N$ is a {\it weak supplement} of a module $L$ in $M$ if
$N+L=M$ and $N \cap L \ll M$ holds.
\end{defn}

We will now state our main result:

\begin{thm} \label{main_theorem} Let $M$ be an $R$-module having
$AB5^{\ast}$ such that every small submodule of $M$ has finite hollow
dimension. Then every submodule of $Rad(M)$ has finite hollow dimension.
\end{thm}

\begin{proof}
Let $G$ be a submodule of $Rad(M)$ with $G \not \ll M$ and assume
$\{ N_\lambda \}_\Lambda$ to be a coindependent family of
proper submodules of $G$ that can be
assumed to be completely coindependent by Lemma \ref{lemma1}. Moreover
we assume that $|\Lambda| \geq 2$.
For all $\lambda \in \Lambda$ there exist elements
$x_\lambda \in Rad(M) \setminus N_\lambda$ since the $N_\lambda$'s are
proper submodules of $G$. Hence $Rx_\lambda \ll M$
and
$L_\lambda := N_\lambda + Rx_\lambda \neq N_\lambda.$
Let $N:=\bigcap_\Lambda N_\lambda$ and $L:=\bigcap_\Lambda L_\lambda$.
Applying Lemma \ref{lemma3}, we get that
$\{(N_\lambda \cap L)/N\}_\Lambda$ is a completely coindependent family of
proper submodules of $L/N$. Next we will show that $L/N$ has finite hollow dimension so that
$\Lambda$ has to be finite.
Since $L_\lambda$ lies above $N_\lambda$ for all $\lambda \in
\Lambda$, we get by applying Lemma \ref{lemma4} that
$L$ lies above $N$ in $M$. Let $K$ be a weak supplement of $L$ in $M$.
Then $L/N \simeq (L\cap K)/(N\cap K)$ yields $hdim(L/N) \leq hdim(L\cap K)$.
By assumption $L\cap K \ll M$ has finite hollow dimension.
Thus $L/N$ has finite hollow dimension and $\Lambda$ must be finite.
This shows that every coindependent family of submodules of $G$ must
be finite. Hence every submodule of $Rad(M)$ has finite hollow dimension.
\end{proof}

Let us recall a result by Lemonnier to prove the next theorem.

\begin{prop} \label{lemon} (Lemonnier, \cite[Proposition
1.3]{lemonnier})\\
Let $M$ be an $R$-module such that every non-zero factor module of $M$
has finite uniform dimension and contains a non-zero submodule having
Krull dimension. Then $M$ has Krull dimension.
\end{prop}

\begin{thm} \label{cor_1} Let $M$ be an $R$-module having $AB5^{\ast}$ such that
every small submodule of $M$ has Krull dimension. Then $Rad(M)$ has Krull dimension.
\end{thm}

\begin{proof}
It is well-known that a module having Krull dimension has finite uniform
dimension (see \cite[6.2]{dung}). Hence every factor module of a small
submodule $N$ of $M$ has finite uniform dimension.
Since $N$ has $AB5^{\ast}$ every submodule of $N$ has finite hollow dimension,
 by Lemma \ref{hs}. Hence
by Theorem \ref{main_theorem} every submodule of $Rad(M)$ has
finite hollow dimension.
By Lemma \ref{hs} every factor module of $Rad(M)$ has finite uniform
dimension. In order to apply Lemonnier's Proposition, we need to show, that
every non-zero factor module of $Rad(M)$ contains a non-zero submodule having Krull
dimension. Let $L\subset Rad(M)$ and $x \in Rad(M)\setminus L$; then $Rx
\ll M$ so that $Rx$ has Krull dimension and hence $(Rx + L)/L \subseteq Rad(M)/L$
has Krull dimension. Applying Proposition \ref{lemon},
$Rad(M)$ has Krull dimension.
\end{proof}

\begin{cor}
Let $M$ be an $R$-module such that $Rad(M)$ has $AB5^{\ast}$ and
every small submodule of $M$ has Krull dimension. Then every submodule of
$Rad(M)$ that has a weak supplement in $M$ has Krull dimension.
\end{cor}

\begin{proof}
 By Theorem \ref{cor_1}, the radical of every submodule contained in
 $Rad(M)$ has Krull dimension. Since $Rad(N)=N \cap Rad(M)$ holds for
 every supplement $N$ in $M$ (see \cite[41.1]{wisbauer}), every supplement
 in $M$ that is a submodule of $Rad(M)$ has Krull dimension.
 Let $L\subseteq Rad(M)$ and $K\subseteq M$ a weak supplement of $L$ in
 $M$. Then $Rad(M)=L + (Rad(M) \cap K)$. Since
 $Rad(M)$ has $AB5^{\ast}$ it is amply supplemented. Thus there exists a
 supplement $N \subseteq L$ of $K\cap Rad(M)$ in $Rad(M)$
 such that $Rad(M)=N + (Rad(M)\cap K)$ and
 $N \cap Rad(M) \cap K = N \cap K \ll N$ holds.
 Moreover $L = N + (L\cap K)$ and $M=N+ K$ holds. Thus $N$ is a
 supplement of $K$ in $M$, implying that $N$ has Krull dimension.
 Because $L/N \simeq (L\cap K)/(N\cap K)$ with $L\cap K \ll M$, $L/N$
 has Krull dimension and hence so too has $L$.
\end{proof}

\begin{center}
ACKNOWLEDGMENTS
\end{center}
The research for this paper was done at the
University of Glasgow. The material in this paper will form part of
the author's MSc thesis for the University of Glasgow and part of
the author's Diplomarbeit for the Heinrich-Heine Universit{\"a}t
D{\"u}sseldorf. The author would like to express his gratitude to
his supervisors: Professor P.F. Smith and Professor R. Wisbauer,
for all their help, advice and encouragement. He would like to thank
especially P.F. Smith for bringing his attention to this problem and for
 helpful discussions about it.

\end{document}